\documentclass[reqno]{amsart}
\usepackage{amsaddr}

\usepackage{amsmath}
\usepackage{amssymb}
\usepackage{amsfonts}
\usepackage{amsthm}
\usepackage{multirow}
\usepackage{enumerate}
\usepackage{color}
\usepackage{dsfont}

\usepackage{setspace}
\newcommand{\R}{\mathbb{R}}
\newcommand{\C}{\mathbb{C}}
\newcommand{\f}{\rightarrow}
\newcommand{\deb}{\ov\partial}
\newcommand{\de}{\partial}
\newcommand{\K}{K\"{a}hler}

\newcommand{\lmb}{\lambda}
\newcommand{\ov}[1]{\overline{#1}}
\newcommand{\w}[1]{\widetilde{#1}}

\newcommand{\D}{\mathcal{D}}

\newcommand{\Aut}{\operatorname{Aut}}

\newcommand{\Ric}{\operatorname{Ric}}

\newcommand{\Ent}{\operatorname{Ent}}

\newcommand{\Isom}{\operatorname{Isom}}

\newcommand{\ol}{\operatorname{Hol}}

\newcommand{\ep}{{\varepsilon}}

\newcommand{\hilb}{\mathcal{H}}

\newcommand{\fdim}{\hspace*{\fill}$\Box$}

\newtheorem{thm}{Theorem}[section]

\newtheorem{lem}[thm]{Lemma}

\newtheorem*{defni}{Definition}

\newtheorem{rmk}[thm]{Remark}

\begin{document}

\author[A. Loi, R. Mossa]{Andrea Loi,  Roberto Mossa}
\address{Dipartimento di Matematica e Informatica, Universit\`{a} di Cagliari,
Via Ospedale 72, 09124 Cagliari, Italy}
\email{loi@unica.it;   roberto.mossa@gmail.com }
\thanks{The authors were  supported by Prin 2010/11 -- Variet\`a reali e complesse: geometria, topologia e analisi armonica -- Italy and also by INdAM. GNSAGA - Gruppo Nazionale per le Strutture Algebriche, Geometriche e le loro Applicazioni.}
\date{\today}
.

\begin{abstract}
In this paper we provide a positive answer to a conjecture due to {A. J. Di Scala}, { A. Loi}, {H. Hishi} (see \cite[Conjecture 1]{dlh}) claiming that a simply-connected homogeneous \K\ manifold $M$ endowed with an integral \K\ form  $\mu_0\omega$, admits a holomorphic isometric immersion in the complex projective space, for a suitable $\mu_0>0$. This result has two corollaries which extend  to   homogeneous K\"ahler manifolds  the results obtained by the authors in \cite{LM03} and in \cite{M04} for homogeneous bounded domains.
\end{abstract}

\title[Some remarks on homogeneous K\"ahler manifolds]
{Some remarks on homogeneous K\"ahler manifolds}


\maketitle

\section{Introduction and statement of the main results}
The paper consists of three results: Theorem \ref{THM04}, Theorem \ref{THM03}  and Theorem   \ref{THM05}. Our first result answer positively to a conjecture due to {A. J. Di Scala}, { A. Loi}, {H. Hishi} (see \cite[Conjecture 1]{dlh}), namely:
\begin{thm}\label{THM04}
Let $(M, \omega)$ be a simply-connected homogeneous \K\ manifold such that its  associated \K\ form  $\omega$ is integral.
Then  there exists a constant $\mu_0  > 0$ such that $\mu_0 \, \omega$  is  projectively induced.
\end{thm}

Before stating our second result (Theorem \ref{THM03}), we recall the definitions of a Berezin quantization and diastasis function. 
Let $(M, \omega)$ be a symplectic manifold and let $\lbrace \cdot, \cdot \rbrace$ be the associated Poisson bracket.
A {\em Berezin quantization} on $M$ is given by a family of associative algebras $\mathcal A_\hbar \subset \C^{\infty}(M)$ such that
\begin{itemize}
\item $\hbar \in E \subset \R^+$ and $\inf E =0$
\item exists a subalgebra ${\mathcal A}\subset (\oplus_{h \in E} \mathcal A_h,*)$, such that for an arbitrary element 
$f=f(\hbar)\in {\mathcal A}$, where $f(\hbar)\in {\mathcal A}_\hbar$, there exists a limit $\lim_{\hbar\rightarrow 0} f(\hbar)=\varphi (f)\in C^{\infty}(M)$,
\item hold the following correspondence principle: for $f, g \in {\mathcal A}$
$$\varphi (f*g)=\varphi (f)\varphi (g),\    \    \  \varphi\left(\hbar^{-1}(f*g-g*f)\right)=i\lbrace\varphi(f), \varphi (g)\rbrace,$$
\item for any pair of points $x_1,x_2 \in \Omega$ there exists $f \in \mathcal A $ such that $\varphi (f) (x_1) \neq \varphi(f)(x_2)$
\end{itemize}

Let $(M,\omega)$ be a real analytic \K\ manifold. Let $U \subset M$ be a neighborhood of a point $p \in M$ and let $\psi: U \f \R$ be a \K\ potential for $\omega$. The potential $\psi$ can be analytically extended to a sesquilinear function $ \psi(p, \ov q)$, defined in a neighborhood of the diagonal of $U \times U$, such that $ \psi(q, \ov q)=\psi(q)$. Assume   (by shrinking $U$)  that the analytic extension $ \psi$ is defined on $U \times U$.  The {\it Calabi's diastasis function} $\D:U \times U \f \R$ is given by:
\[
\D(p,q)=\psi(p,\ov p)+\psi(q,\ov q)+\psi(p,\ov q)+\psi(q,\ov p).
\] 
Denoted by $\D_p(q):=\D(p,q)$ the diastasis centered in a point $p$, one can see that $\D_p:U \f \R$ is a \K\ potential around $p$.

Our second result extends to any homogeneous \K\ manifolds the results obtained by the authors  \cite{LM03} for homogeneous bounded domains.

\begin{thm}\label{THM03}
Let $(M, \omega)$ be a homogeneous \K\ manifold. Then the following are equivalent:
\begin{enumerate}[(a)]
\item $M$ is contractible.
\item $(M, \omega)$ admits a global \K\ potential.
\item $(M, \omega)$ admits a global diastasis $\D: M \times M \f \R$.
\item $(M, \omega)$ admits a Berezin quantization.
\end{enumerate} 
\end{thm}

Also to state   the third result (Theorem \ref{THM05}) we need some definitions.
The diastatic entropy has been defined by the second author in \cite{M04} (see also \cite{mossa1}) following the ideas developed in \cite{LM01} and \cite{mossa2}.
Assume that $\omega$ admits a globally defined  diastasis function  $\D_p:M \f \R$ (centered at $p$). The \emph{diastatic entropy} at $p$ is defined as
\begin{equation}\label{def dent}
\Ent\left( M, \, \omega \right)\left( p\right)= \min \left\{ c >0 \ | \, \int_M e^{-c \, \D_p }\, \frac{\omega^n}{n!}<\infty \right\} .
\end{equation}
The definition does not depend on the point $p$ chosen  (see \cite[Proposition 2.2]{mossa1}).

Assume that $M$ is simply-connected and assume that there exists a line bundle $L \f M$ such that $c_1(L)=[\omega]$ (i.e. $\omega$ is integral). Let $h$ be an Hermitian metric on $L$ such that $\Ric(h)=\omega$ and consider the Hilbert space of global holomorphic sections of $L^{\lmb}=\otimes^\lmb L$ given by
\begin{equation}\label{hilbertspaceh}
\hilb_{\lmb,h}=\left\{ s\in\ol(L) \ |\   \|  s \|^2 = \int_M  h^{\lmb}(s,s)\, \frac{\omega^n}{n!}<\infty\right\},
\end{equation}
with the scalar product induced by $\| \cdot \|$.
Let $\{s_j\}_{j=0,\dots,N}$, $N \leq \infty$, be an orthonormal basis for $\hilb_{\lmb, h}$. The \emph{$\varepsilon$-function} is given by
\begin{equation}\label{epl}
\varepsilon_{\lmb}(x)=\sum_{j=0}^Nh^\lmb(s_j(x),s_j(x)).
\end{equation}
This definition depends only on the \K\ form $\omega$. Indeed since $M$ is simply-connected, there exists (up to isomorphism) a unique $L \f M$ such that $c_1(L)=[\omega]$, and 
 it is easy to see that the definition  does not depend on the orthonormal basis chosen or on the Hermitian metric $h$ (see e.g \cite{LM03} or \cite{M01,LM02} for details). 
Under the assumption that $\varepsilon_{\lmb}$ is a (strictly) positive function, one can  consider the  \emph{coherent states map} $f: M \f \C P^N$ defined by
\begin{equation}\label{eqcoherent}
f(x)=\left[s_0(x), \dots, s_j(x), \dots\right].
\end{equation}
The fundamental  link between the coherent states map and the $\varepsilon$-function
is expressed by the following equation (see \cite{rawnsley} for a proof)
\begin{equation}\label{pullbackdie}
f^*\omega_{FS}=\lmb \omega + \frac{i}{2} \de \deb\, \varepsilon_\lmb ,
\end{equation}
where $\omega_{FS}$ is the Fubini--Study form on $\C P^N$.
 
\begin{defni}
 We say that $\lambda\omega$ is a  \emph{balanced metric} if and only if the $\varepsilon_\lmb$ is a \emph{positive constant}.
 \end{defni}

We can now state our  third and last result which  extends to any homogeneous \K\ manifold, the result obtained by the second author \cite[Theorem 2]{M04}, about the link between \emph{diastatic entropy} and \emph{balanced condition} on homogeneous bounded domains.

\begin{thm}\label{THM05}
Let $(M,\omega)$ be a contractible homogeneous \K\ manifold. Then $\lmb \omega$ is balanced if and only if 
\begin{equation}
\Ent(M,\omega)<\lmb.
\end{equation}
\end{thm}

\section{Proof the main results}
We start with the following two lemmata, which provide a necessary and sufficient condition on the non-triviality of the Hilbert space $\hilb_{\lmb, h}$.\begin{lem}[Rosemberg--Vergne \cite{rosenberg}]\label{lemrose}\footnote{The authors thanks Hishi Hideyuki for reporting this result.}
Let $(M, \omega)$ be a simply-connected homogeneous \K\ manifold with $\omega$ integral. Then there exists $\lambda_0>0$ such that  
$\hilb_{\lmb, h}\neq \{0\}$ if and only if $\lambda >  \lambda_0$ and $\lmb \omega$ is integral.
\end{lem}

\begin{lem}   \label{lembal}
Let $(M, \omega)$ be a  simply-connected homogeneous \K\ manifold. Then   $\hilb_{\lambda,h} \neq \{0\}$ if and only if $\lambda \omega$ is a balanced metric.
\end{lem}
\proof Let $F \in \Aut(M)\cap\Isom(M,\omega)$ be an holomorphic isometry and let $\w F$ its lift to $L$ (which exists since $M$ is simply-connected). Note that, if $\lbrace s_0, \dots,s_N\rbrace$, $N\leq\infty$, is an orthonormal basis for $\hilb_{\lmb, h}$, then $\lbrace\w{ F}^{-1} \left( s_0\left(F\left(x\right) \right)\right),\dots , \w{ F}^{-1} \left( s_N\left(F\left(x\right) \right)\right)\rbrace$ is an orthonormal basis for $\hilb_{\lmb,\w{ F} ^*h}$. 
Therefore
\begin{equation*}
\begin{split}
\epsilon_\lmb(x)
&=\sum_{j=0}^N \w {F}^* h^\lmb \left( \w {F}^{-1}(s_j(F(x))),\w {F}^{-1}(s_j(F(x)))\right) \\
&=\sum_{j=0}^N h^\lmb\left( s_j(F(x)), s_j(F(x)) \right)=\varepsilon_{\lmb}\left(F(x)\right).
\end{split}
\end{equation*}
Since $\Aut(M)\cap\Isom(M,\omega)$ acts transitively on $M$, $\varepsilon_{\lmb}$ is forced to be constant. 
\endproof

\noindent
{\bf Proof of Theorem \ref{THM04}}\label{proofthm1}
By Lemma \ref{lemrose} 
there exists $\lmb > \lmb_0$, such that the Hilbert space 
$\hilb_{\lmb, h}\neq \{0\}$.
By Lemma \ref{lembal}, $\varepsilon_\lmb$ is a positive constant, so the coherent states map $f$ given by \eqref{eqcoherent} is well defined. Moreover, by \eqref{pullbackdie}, we have that $f^*\omega_{FS}=\lmb \omega$, i.e. $\lmb \omega$ is projectively induced  for all $\lmb > \lmb_0$. The conclusion follows by setting $\mu_0 >\lmb_0$. \fdim

\vskip 0.5cm

The main ingredients for the proof of Theorem \ref{THM03} are the following two lemmata.
The first one is the celebrated theorem of Dorfmeister and Nakajima which provides a positive answer to the so called \emph{Fundamental Conjecture} formulated by Vinberg and Gindikin. The second one  is a reformulation due to  Engli\v{s}  of the Berezin   quantization result       in terms of the $\ep$-function   and  Calabi's diastasis function.

\begin{lem}[Dorfmeister--Nakajima \cite{DN88}]\label{THMA}
A homogeneous \K\ manifold $(M, \omega)$  is the total space of a holomorphic
fiber bundle over a homogeneous bounded domain $\Omega$ in which the fiber ${\mathcal F} ={\mathcal E} \times {\mathcal C}$
is (with the
induced K\"ahler metric) is the  \K\ product of a flat homogeneous K\"ahler manifold ${\mathcal E}$ and a
compact simply-connected homogeneous K\"ahler manifold ${\mathcal C}$.
\end{lem}

\begin{lem}[Engli\v{s} \cite{englis}]\label{thmberezin}
Let $\Omega\subset \C^n$ be a complex domain equipped with a real analytic  \K\ form $\omega$.
Then, $(\Omega, \omega)$ admits a Berezin quantization if the following two  conditions are satisfied:
\begin{enumerate}
\item
the function $\ep_{\lambda}(x)$ is a positive constant (i.e. $\lmb \omega$ is balanced) for all sufficiently large $\lambda$;
\item
the function $e^{-\D(x, y)}$  is globally defined on $\Omega\times\Omega$,  $e^{-\D(x, y)}\leq1$  and $e^{-\D(x, y)}=1$,  if and only if $x=y$.
\end{enumerate}
\end{lem}

\noindent
{\bf Proof of Theorem \ref{THM03}} 
\noindent
\emph{\bf (a) $\Rightarrow$ (b)}.
By Lemma \ref{THMA}, since a homogeneous bounded domain is contractible, $M$ is a {\em complex}
product $\Omega \times {\mathcal F}$, where  ${\mathcal F}={\mathcal E}\times {\mathcal C}$ is (with the induced metric) a \K\ product  of a flat  \K\ manifold ${\mathcal E}$ and a compact simply-connected homogeneous \K\ manifold ${\mathcal C}$. On the other hand
${\mathcal E}$ is \K\ flat and therefore ${\mathcal E}=\C^k\times T$ where $T$  is a product of flat complex tori. Hence $M$ is a {complex}
product $\Omega \times \C^k\times T \times {\mathcal C}$. Since we assumed $M$ contractible, the compact factor $T \times {\mathcal C}$ has dimension zero and $M=\Omega \times \C^k$.  It is well-know that $\Omega$ is biholomorphic to a Siegel domain  (see, \cite{vinberg} for a proof), therefore  $\Omega \times \C^k$ is pseudoconvex and a classical result of Hormander (see \cite{bulletin}) asserting that the equation $\deb u = f$ with $f$ $\deb$-closed form has a global solution on pseudoconvex domains, assures us the existence of a global potential $\psi$ for $\omega$ (see also \cite{mossa3, M04}, and the proof of Theorem 4 in \cite{dlh} for an explicit construction of the potential $\psi$).

\emph{\bf (b) $\Rightarrow$ (c).} Let $\psi: M \f \R$ be a global \K\ potential for $(M, \omega)$. By Lemma \ref{THMA}, $M=\Omega \times \C^k\times T \times {\mathcal C}$. The compact factor $T \times {\mathcal C}$ is a \K\ submanifold of $M$, therefore the existence of a global \K\ potential on $M$ implies that $\dim(T \times {\mathcal C})=0$. So $M=\Omega\times \C^k$. 

Consider the Hilbert space $\hilb_{\lmb, h}$ defined in \eqref{hilbertspaceh}. 
Since $\Omega\times \C^k$ is contractible the line bundle $L \cong M \times \C$. So, the holomorphic section $s \in \hilb_{\lmb, h}$ can be viewed as a holomorphic function $s:M \f \C$. As Hermitian metric $h$ over $L$ we can take the one defined by $h(s,s)= e^{-\psi} |s|^2$. Hence $\hilb_{\lmb, h}$ can be identified with the weighted Hilbert space $\hilb_{\lmb\psi}$ (see \cite{LMZ03}) of square integrable holomorphic functions on $M=\Omega \times \C^k$, with weight $e^{-\lmb\psi}$, namely
\begin{equation}\label{hilbertspacePhi}
\hilb_{\lmb\psi}=\left\{ s\in\ol(M) \ | \ \, \int_M e^{-\lmb\psi}|s|^2\frac{\omega^n}{n!}<\infty\right\}.
\end{equation}
Assume $\lmb > \lmb_0$ with $\lmb_0$ given by  Lemma \ref{lemrose}, so that $\hilb_{\lmb\psi}\neq \{0\}$. Let $\{s_j\}$ be an orthonormal basis for $\hilb_{\lmb\psi}$, then the reproducing kernel is given by
$$K_{\lmb\psi}(z, \bar w)=\sum_{j=0}^\infty s_j(z) \ov{s_j(w)} .$$
Then, the $\ep$-function (defined in \eqref{epl}) reads as: 
\begin{equation}\label{epsilon}
\varepsilon_\lmb(z)=e^{-\lmb\psi(z)}K_{\lmb\psi}(z,\bar  z).
\end{equation}
Let $ \psi \left( z,\ov w\right) $ be the analytic continuation of the \K\ potential $\psi$. By Lemma \ref{lembal} there exists a constant $C$ such that

\begin{equation}\label{epsilon1}
\ep_{\lmb}\left( z,\, \ov w\right) =e^{-\lmb\psi\left( z,\, \ov w\right) }K_{{\lmb\psi}}\left( z,\, \ov w\right) =C>0.
\end{equation}
Hence $K_{{\lmb\psi}}\left( z,\, \ov w\right)$ never vanish. Therefore, for any fixed point $z_0$, the function
\begin{equation}\label{PSI00}
\Phi\left( z,\, \ov w\right)=\frac{K_{{\lmb\psi}}\left( z,\,\ov w\right) K_{{\lmb\psi}}\left( z_0,\,\ov z_0\right)  }{K_{{\lmb\psi}}\left( z,\,\ov z_0\right) K_{{\lmb\psi}}\left( z_0,\,\ov w\right) }
\end{equation}
is well defined. Note that 
\begin{equation}\label{DIAST00}
\D_{z_0}\left( z\right) =\frac{1}{\lmb}\log\Phi\left( z,\,\ov z\right) 
\end{equation}
is  the diastasis centered in $z_0$ associated to $\omega$ and that $\D_{z_0}$ is defined on whole $M$. Since we can repeat this argument for any $z_0 \in M$, we conclude that the diastasis $\D : M \times M \f \R$ is globally defined.

\emph{\bf (c) $\Rightarrow$ (d).} Arguing as in \lq\lq{\bf (b) $\Rightarrow$ (c)}\rq\rq, the existence of a global diastasis implies that $M$ is a complex product $\Omega \times \C^k$. Therefore, as in \eqref{hilbertspacePhi} $\hilb_{\lmb, h}\cong\hilb_{\lmb \D_{z_0}}=\left\{ s\in\ol(M) \ | \ \, \int_M e^{-\lmb\\D_{z_0}x}|s|^2\frac{\omega^n}{n!}<\infty\right\}$. Assume $\lmb > \lmb_0$ with $\lmb_0$ given by  Lemma \ref{lemrose} and consider the {coherent states map} $f$ given by \eqref{eqcoherent}. By Lemma \ref{lembal}  $\varepsilon_\lmb$ is a positive constant and by \eqref{pullbackdie} we conclude that $f^*\omega_{FS}= \lmb \omega$.

By \cite[Example 6]{LM03}, the
Calabi's diastasis function $\D_{{FS}}$ associated to  $\omega_{FS}$ is such that  $e^{-\D_{FS}}$
is globally defined on $\C P^{N}\times\C P^{N}$. Since the diastasis $\D$ is globally defined on $M$, by the hereditary property of the diastasis function
 (see \cite[Proposition  6]{calabi}) we get that, for all $x, y\in M$,
\begin{equation}\label{fondequ}
e^{-\D_{FS}(f(x), f(y) )}= e^{-\lambda \D(x, y)}=\left(e^{- \D(x, y)}\right)^{\lambda}
\end{equation}
 is globally defined on $M\times M$.
 Since, by \cite[Example 6]{LM03},   $e^{-\D_{FS}(p, q)}\leq 1$ for all $p, q\in \C P^{N}$ it follows that $e^{-\D(x, y)}\leq 1$ for  all
 $x, y\in M$.
By Lemma \ref{thmberezin}, it remains to show that   $e^{-\D(x, y)}=1$ iff $x=y$.  By (\ref{fondequ})   
and by the fact that $e^{-\D_{FS}(p, q)}=1$ iff $p=q$ (again by \cite[Example 6]{LM03})
this is equivalent to   the injectivity of the coherent states map  $f$. This follows by \cite[Theorem 3]{dlh}, which asserts that   a   \K\ immersion  of a homogeneous \K\ manifold into a  finite or infinite dimensional complex projective space is  one to one.

\emph{\bf (d) $\Rightarrow$ (a).} By the very definition of Berezin quantization there exists a global potential for $(M, \omega)$. By Lemma \ref{THMA} we deduce, as above, that $M$ is a complex product $\Omega \times \C^k$, where $\Omega$ is a bounded homogeneous domain which is contractible.\fdim

\vskip 0.5cm

\noindent
{\bf Proof of Theorem \ref{THM05}}
\noindent
By (c) in Theorem \ref{THM03}, the diastasis $\D: M \times M \f \R$ is globally defined. Assume that $\lmb\, \omega$ is balanced i.e. that $\ep_\lmb$ is a positive constant. Since the $\ep_\lmb$ does not depend on the \K\ potential, by \eqref{epsilon1} we have
\begin{equation}\label{epsilon2}
\ep_{\lmb}\left( z,\, \ov w\right)=e^{-\lmb\,\D_{z_0}\left( z,\, \ov w\right)}K_{\lmb\,\D_{z_0}}\left( z,\, \ov w\right)=C
\end{equation}
where $\D_{z_0}\left( z,\, \ov w\right)$ denote the analytic continuation of $\D_{z_0}\left( z\right)$ and $K_{\lmb\D_{z_0}}$ is the reproducing kernel for $\hilb_{\lmb\D_{z_0}}$. By \eqref{PSI00}, with $\ov w = \ov z_0$, we get $\D_{z_0}\left( z,\, \ov z_0\right)= \frac{1}{\lmb} \log\Phi\left( z,\,\ov z_0\right)=0$. Hence
\begin{equation*}
C= e^{-\lmb\,\D_{z_0}\left( z,\, \ov z_0\right)} K_{\lmb\,\D_{z_0}}\left( z,\, \ov z_0\right)=K_{\lmb\,\D_{z_0}}\left( z,\, \ov  z_0\right)  \in \hilb_{\lmb\,\D_{z_0}}.
\end{equation*}
Thus $\hilb_{{\lmb\,\D_{z_0}}}$ contains the constant functions and
\[
\int_M e^{-\lmb\,\D_{z_0}}\ \frac{\omega^n}{n!}<\infty.
\]
Therefore, by definition of diastatic entropy, 
\[
\Ent \left( M ,\omega\right)\left( z_0 \right) <\lmb<\infty.
\]
On the other hand, if for some $z_0$, $\Ent \left( M,\omega\right)(z_0)<\lmb$, then $\hilb_{\lmb\D_{z_0}}\not =\lbrace 0 \rbrace$ and by Lemma \ref{lembal} we conclude that $\lmb \omega$ is balanced. \fdim

\begin{rmk}\rm
From the previous proof, we see that if $M$ is simply-connected, then $\Ent(M, \omega)(z_0) = \lmb_0$ for any $z_0 \in M$, where $\lmb_0$ is the positive constant defined in Lemma \ref{lemrose}. 
\end{rmk}

\end{document}